\documentclass{article}
\usepackage{amsmath}
\usepackage{amssymb}
\usepackage{amsthm} 
\usepackage{graphicx}
\usepackage{multirow}
\usepackage{cite}
\usepackage{enumerate}
\usepackage[OT2,T1]{fontenc}
\DeclareSymbolFont{cyrletters}{OT2}{wncyr}{m}{n}
\DeclareMathSymbol{\Sha}{\mathalpha}{cyrletters}{"58}

\usepackage{amsfonts}

\newcommand{\Z}{\mathbb{Z}}
\newcommand{\Q}{\mathbb{Q}}

\newcommand{\Qp}{\mathbb{Q}_p}

\begin{document}

\title{On the Rank of the Elliptic Curve $y^2=x(x-p)(x-2)$}  
\author{Jeffrey Hatley\footnote{Department of Mathematics and Statistics, University of Massachusetts - Amherst}}
\maketitle

\begin{abstract}
An elliptic curve $E$ defined over $\Q$ is an algebraic variety which forms a finitely generated abelian group, and the structure theorem then implies that $E \cong \Z^{r} \oplus \Z_{\text{tors}}$ for some $r \geq 0$; this value $r$ is called the rank of $E$. It is a classical problem in the study of elliptic curves to classify curves by their rank. In this paper, the author uses the method of $2-$descent to calculate the rank of two families of elliptic curves.
\end{abstract}

\section{Introduction}

An elliptic curve $E$ defined over $\Q$ is the set of solutions $(x,y) \in \Q^2$ to an equation of the form
\[
y^2 = x^3 + ax^2 + bx + c, \text{ with } a,b,c \in \Q
\]
\noindent along with an additional point at infinity, ${\cal{O}}$. A classical theorem from Mordell shows that such an ellipic curve forms a finitely generated abelian group; the structure theorem then implies that 
\[
E \cong \Z^{r} \oplus \Z_{\text{tors}}, \text{   } r \geq 0.
\] 
\noindent This nonnegative integer $r$ is called the \textit{rank} of the elliptic curve. It is a classical problem to classify ellipic curves, and the rank of a curve provides a useful way to distinguish it from other curves, as well as to gain some insight into its algebraic structure. However, calculating the rank of a given ellipic curve can be quite difficult in general. One method of doing so is the method of 2-descent. The machinery behind 2-descent is quite high-powered, but the idea is rather simple. 

The main idea behind $2$-descent is a local-to-global method. One indirectly studies an elliptic curve $E$ by examining whether certain related equations, called \textit{homogeneous spaces}, have rational points over every local field $\Qp$ for $p = \infty$ or $p$ a rational prime. This information is then pieced together to yield information about $E$. 

This paper uses two variations on the method of two-descent to gain results on the ranks of two different families of ellipic curves. While studying computer computations of the ranks of elliptic curves of the form
\[
E_p : y^2 = x(x-p)(x-2),
\]
Jason Beers made the following conjecture.

\bigskip

\textbf{Conjecture. } \textit{Let $E_p$ be the elliptic curve defined by 
\\ $E_p : y^2 = x(x-p)(x-2)$ where $p$ and $p-2$ are twin primes. Then}
\begin{equation*}
\text{rank}(E_p) = 
\begin{cases} 0 & \text{if $p \equiv 7$ (mod $8$),}
\\
1 & \text{if $p \equiv 3,5$ (mod $8$),}
\\
2 & \text{if $p \equiv 1$ (mod $8$).}
\end{cases}
\end{equation*}

Curves similar to this form were considered in Silverman [1], and similar results for the curve $E:y^2= x(x-(p-2))(x-p)$ were recently made available on arXiv.

\bigskip

In this paper, we prove the first case and part of the second; more specifically, we prove

\bigskip

\textbf{Theorem 1. } \textit{Let $p$ and $p-2$ be twin primes, and let $E_p$ be the elliptic curve defined by $E_p:y^2=x(x-p)(x-2)$. }

(a) \textit{If $p \equiv 7$ (mod $8$), then rank$(E_p)=0.$ In particular, we have}
\[
E_p(\Q) \cong \Z / 2\Z \times \Z / 2\Z.
\]

(b) \textit{If $p \equiv 5$ (mod $8$), then rank$(E_p) \leq 1$. In particular, we have}
\[
E_p(\Q) \cong \Z / 2\Z \times \Z / 2\Z.
 \]
\text{or}
\[
E_p(\Q) \cong \Z \times \Z / 2\Z \times \Z / 2\Z.
 \]

\bigskip

Before going through each proof, we state the method of two-descent to be used as presented in Silverman [1].

\section{Proof of Theorem 1(a)}

In this section we prove the following theorem:

\bigskip

\textbf{Theorem 1(a).} \textit{Let $p$ and $p-2$ be prime numbers in $\Z$ with $p \equiv 7$ (mod $8$). Then the elliptic curve $E(\Q)$ given by}
\[
E: y^2 = x(x-p)(x-2)
\] 
\textit{has rank $0$. In particular, }
\[
E(\Q) \cong \Z / 2\Z \times \Z / 2\Z.
 \]
 
\bigskip

To prove this theorem, we use the method of $2$-descent presented as Proposition 1.4 of Chapter X in Silverman [1], which we now state. 

\bigskip

\textbf{Theorem 2. (Complete 2-Descent, Version 1).} \textit{Let $E/\Q$ be an elliptic curve given by an equation} 
\[
y^2 = (x-e_1)(x-e_2)(x-e_3) \text{     with } e_{1},e_{2},e_{3} \in \Q.
\]
\textit{Let $S$ be a set of places of $\Q$ including $2, \infty$, and all places dividing the discriminant of $E$. Further, let }
\[
\Q(S,2) = \{ b \in \Q^* / \Q^{*2} : \text{ord}_{\nu}(b) \equiv 0 \text{(mod $2$) for all } \nu \notin S\}.
\]
\textit{There is an injective homomorphism}
\[
E(\Q) / 2E(\Q) \to \Q(S,2) \times \Q(S,2)
\]
\textit{defined by}
\[
P=(x,y) \to \begin{cases}

  (x-e_1,x-e_2) & \textit{if $x \neq e_{q},e_{2}$} \\

  ((e_{1}-e_{3})/(e_{1}-e_{2}),e_{1}-e_{2}) & \textit{if $x=e_1$} \\

  (e_2 - e_1, (e_2 - e_3)/(e_2 - e_1)) & \textit{if $x=e_2$} \\
  
  (1,1) & \textit{if $x= \infty$ (i.e. if $P={\cal{O}}$}).

\end{cases}
\]
\textit{Let $(b_1,b_2) \in \Q(S,2) \times \Q(S,2)$ be a pair which is not the image of one of the three points ${\cal{O}}, (e_1, 0), (e_2,0)$. Then $(b_1,b_2)$ is the image of a point \\ $P=(x,y) \in E(\Q) / 2E(\Q)$ if and only if the equations }
\[
b_{1}z_{1}^2 - b_{2}z_{2}^2 = e_{2} - e_{1}
\]
\[
b_{1}z_{1}^2 - b_{1}b_{2}z_{3}^2 = e_{3} - e_{1}
\]
\textit{have a solution $(z_{1},z_{2},z_{3}) \in \Q^* \times \Q^* \times \Q^*$; if such a solution exists, then one can take }
\[
P = (x,y) = (b_{1}z_{1}^2 + e_{1}, b_{1}b_{2}z_{1}z_{2}z_{3}).
\]

\bigskip

Thus, Theorem $2$ allows one to calculate $E(\Q) / 2E(\Q)$ for ellipic curves defined by a sufficiently nice equation. If $E(\Q) \cong \Z^{r} \oplus \Z_{\text{tors}}$, then all of the odd-torsion (i.e. factors of $E(\Q)$ of the form $\Z / m\Z$ where $m$ is odd) is killed in $E(\Q) / 2E(\Q)$. Furthermore, if $p$ is a prime of good reduction for $E$ and gcd($p,m$)=$1$, then $E(\Q)[m]$ injects into $\tilde{E}(\mathbb{F}_{p})$, the $\mathbb{F}_p$ rational points on the reduction of $E$ modulo $p$, so it is easy to calculate $E_{\text{tors}}$. These facts together with Theorem $2$ allow one to calculate the rank of certain elliptic curves. We now prove Theorem 1(a).

\bigskip

\bigskip

\textbf{Proof of Theorem 1(a).}

Our curve is 
\[
E: y^2 = x(x-p)(x-2) = x^3 -(p+2)x^2 + 2px
\]
which has discriminant $\Delta = 2^{2}p^{2}(p-2)^{2}$, so our set $S$ is
\[
S =\{ 2,p,p-2\}
\]
and our set $\Q(S,2)$ is
\[
\Q(S,2) = \{ \pm 1, \pm 2, \pm p, \pm (p-2), \pm 2p, \pm 2(p-2), \pm p(p-2), \pm 2p(p-2)\}.
\]
The two-torsion points of $E(\Q)$ are those points $(x,y) \in E(\Q)$ with $y=0$, hence the two-torsion of our curve is 
\[
E(\Q)[2] = \{{\cal{O}}, (0,0), (2,0), (p,0) \}.
\]
Now, since $p \geq 7$, we have $3 \nmid \Delta$, and so we see that $E_{\text{tors}}$ injects into $E(\mathbb{F}_3)$. By hypothesis $p=7 + 8k$ for some non-negative integer $k$. In fact, $k \equiv 0$ (mod $3$), since if $k \equiv  1$ (mod $3$) then $p \equiv 0$ (mod $3$) and so is not prime; likewise, if $k \equiv 2$ (mod $3$), then $p-2 \equiv 0$ (mod $3$), so $p-2$ is not prime. Hence we have $k \equiv 0$ (mod $3$), and 
\[
E(\mathbb{F}_{3}) = \{ {\cal{O}}, (0,0), (1,0), (2,0) \}
\]
where the fact that $(1,0) \in E$ follows from our congruence argument on $k$. Any odd $m-$torsion must be an $m-$group, and since $E[2] \subset E_{\text{tors}}$ and both $E[2]$ and $E(\mathbb{F}_3)$ have cardinality $4$, we see that $E_{\text{tors}}=E[2]$.

Now we consider the map 
\[
\phi: E(\Q)/2E(\Q) \to \Q(S,2) \times \Q(S,2)
\]
given in Theorem $2$ with $e_{1}=0$, $e_{2}=2$, and $e_{3}=p$. There are $256$ pairs $(b_1,b_2) \in \Q(S,2) \times \Q(S,2)$, and for each pair we must check to see whether it comes from an element of $\Q(S,2) \times \Q(S,2)$. Using Theorem 2, we can compute the image  $\phi(E[2])$ in $\Q(S,2) \times \Q(S,2)$:
\[
{\cal{O}} \to (1,1) \quad (0,0) \to (2p,-2) \quad (2,0) \to (2,-2(p-2)) \quad (p,0) \to (p,p-2).
\] 
For the remaining pairs $(b_1,b_2)$ we must determine whether the equations
\begin{equation}\label{eq1}
b_{1}z_{1}^{2} - b_{2}z_{2}^{2} = 2
\end{equation}
and
\begin{equation}\label{eq2}
b_{1}z_{1}^{2} - b_{1}b_{2}z_{3}^2 =p
\end{equation}
have a simultaneous solution $(z_{1},z_{2},z_{3}) \in \Q^3$. This is facilitated by a few facts. First, recall that $\Q \subset \Q_{q}$ (the $q$-adic completion of $\Q$) for each prime $q$, so if an equation has no solutions over $\Q_{q}$ then it has no solutions over $\Q$, and the pair $(b_{1},b_{2})$ in that case is not in $\Q(S,2) \times \Q(S,2)$. Second, the map $\phi$ is a homomorphism. Thus, if $(b_{1},b_{2})$ and $(b_{1}',b_{2}')$ are both in the image of $\phi$, then so is $(b_{1}b_{1}',b_{2}b_{2}')$; if $(b_{1},b_{2})$ is in the image and $(b_{1}',b_{2}')$ is not, then $(b_{1}b_{1}',b_{2}b_{2}')$ is not. If the equations corresponding to a pair $(b_{1},b_{2})$ have no solutions over some $\Q_{p}$, we say that $(b_{1},b_{2})$ is $\Q_{p}$\textit{-non-trivial}.

We follow Silverman's lead and list our results in a table whose entries list either the point $(x,y)\in E(\Q)$ that gets mapped to the pair $(b_{1},b_{2})$ or the field over which the equations $\eqref{eq1}$ and $\eqref{eq2}$ have no solution. If $(z_{1},z_{2},z_{3})$ is a solution to equations $\eqref{eq1}$ and $\eqref{eq2}$, then the pre-image of $(b_{1},b_{2})$ is $(b_{1}z_{1}^2 + e_{1}, b_{1}b_{2}z_{1}z_{2}z_{3})$. Each table entry also has a superscript number $(n)$ which refers to a note explaining the entry. We do exclude half of the points from the table, however: it is easy to see that if $b_{1} < 0$ and $b_{2} > 0$ then equation $\eqref{eq1}$ has no solutions in $\mathbb{R}$, and if $b_{1} <0$ and $b_{2} < 0$ then equation $\eqref{eq2}$ has no solution in $\mathbb{R}$. Hence, we exclude the portion of the table with $b_{1} < 0$.

\bigskip


\bigskip

\begin{center}
\begin{tabular}{|c||c|c|c|c|c|c|c|c|} 
\hline
$b_{1}$ \textbackslash  $b_{2}$ & $1$ & $2$ & $p$ & $p-2$ & \ \ $2p$ \ \ & $2(p-2)$ & $p(p-2)$ & $2p(p-2)$ \\ 
\hline
\hline
$1$ & ${\cal{O}}$ & $\Q_{2}^{(3)}$ & $\Q_{p}^{(5)}$ & $\Q_{p-2}^{(7)}$ & \multicolumn{2}{|c|}{$\Q_{2}^{(3)}$}  & $\Q_{p}^{(5)}$ & $\Q_{2}^{(3)}$ \\
\hline
$2$ & $\Q_{2}^{(4)}$ & $\Q_{p-2}^{(12)}$ & \multicolumn{2}{|c|}{$\Q_{2}^{(4)}$} & $\Q_{p}^{(16)}$ & $\Q_{p-2}^{(16)}$ & $\Q_{2}^{(4)}$ & $\Q_{p}^{(14)}$ \\
\hline
$p$  & $\Q_{p}^{(9)}$ & \multirow{2}{*}{$\Q_{2}^{(3)}$} & $\Q_{p}^{(12)}$ & $\Q_{p}^{(9)}$ & \multicolumn{2}{|c|}{\multirow{2}{*}{$\Q_{2}^{(3)}$}} & $\Q_{p}^{(10)}$ & \multirow{2}{*}{$\Q_{2}^{(3)}$} \\
\cline{1-2}\cline{4-5}\cline{8-8}
$p-2$ & $\Q_{p}^{(6)}$ & & $(p,0)$ & $\Q_{p}^{(6)}$ & \multicolumn{2}{|c|}{} & $\Q_{p-2}^{(3)}$ & \\
\hline
$2p$ & \multirow{2}{*}{$\Q_{2}^{(4)}$} & $\Q_{p}^{(9)}$ & \multicolumn{2}{|c|}{\multirow{2}{*}{$\Q_{2}^{(4)}$}} & $\Q_{p}^{(10)}$ & $\Q_{p}^{(9)}$ & \multirow{2}{*}{$\Q_{2}^{(4)}$} & $\Q_{p}^{(10)}$  \\ 
\cline{1-1}\cline{3-3}\cline{6-7}\cline{9-9}
$2(p-2)$ & & $\Q_{p}^{(14)}$ & \multicolumn{2}{|c|}{} & $\Q_{p-2}^{(12)}$ & $\Q_{p}^{(14)}$ &  & $\Q_{p-2}^{(16)}$ \\
\hline  
$p(p-2)$ & $\Q_{p}^{(9)}$ & $\Q_{2}^{(3)}$ & $\Q_{p}^{(10)}$ & $\Q_{p}^{(9)}$ & \multicolumn{2}{|c|}{$\Q_{2}^{(3)}$} & \multirow{2}{*}{$\Q_{2}^{(4)}$} & $\Q_{2}^{(3)}$ \\
\cline{1-7}\cline{9-9}
$2p(p-2)$ & $\Q_{2}^{(4)}$ & $\Q_{p}^{(9)}$ & \multicolumn{2}{|c|}{$\Q_{2}^{(4)}$} & $\Q_{p}^{(10)}$ & $\Q_{p}^{(9)}$ & & $\Q_{p}^{(10)}$ \\ 
\hline
$-1$ & $\Q_{p}^{(14)}$ & $\Q_{2}^{(3)}$ & $\Q_{p-2}^{(12)}$ & $\Q_{p}^{(14)}$ & \multicolumn{2}{|c|}{$\Q_{2}^{(3)}$}  & $\Q_{p-2}^{(15)}$ & $\Q_{2}^{(3)}$ \\
\hline
$-2$ & $\Q_{2}^{(4)}$ & $\Q_{p}^{(6)}$ & \multicolumn{2}{|c|}{$\Q_{2}^{(4)}$} & $(0,0)$ & $\Q_{p}^{(6)}$ & $\Q_{2}^{(4)}$ & $\Q_{p-2}^{(8)}$ \\
\hline
$-p$  & $\Q_{p}^{(9)}$ & \multirow{2}{*}{$\Q_{2}^{(3)}$} & $\Q_{p}^{(10)}$ & $\Q_{p}^{(11)}$ & \multicolumn{2}{|c|}{\multirow{2}{*}{$\Q_{2}^{(3)}$}} & $\Q_{p}^{(10)}$ & \multirow{2}{*}{$\Q_{2}^{(3)}$} \\
\cline{1-2}\cline{4-5}\cline{8-8}
$-(p-2)$ & $\Q_{p-2}^{(11)}$ & & $\Q_{p}^{(13)}$ & $\Q_{p-2}^{(16)}$ & \multicolumn{2}{|c|}{} & $\Q_{p}^{(13)}$ & \\
\hline
$-2p$ & \multirow{2}{*}{$\Q_{2}^{(4)}$} & $\Q_{p}^{(9)}$ & \multicolumn{2}{|c|}{\multirow{2}{*}{$\Q_{2}^{(4)}$}} & $\Q_{p}^{(10)}$ & $\Q_{p}^{(9)}$ & \multirow{2}{*}{$\Q_{2}^{(4)}$} & $\Q_{p}^{(10)}$  \\ 
\cline{1-1}\cline{3-3}\cline{6-7}\cline{9-9}
$-2(p-2)$ & & $(2,0)$ & \multicolumn{2}{|c|}{} & $\Q_{p}^{(6)}$ & $\Q_{p-2}^{(8)}$ &  & $\Q_{p}^{(6)}$ \\
\hline  
$-p(p-2)$ & $\Q_{p}^{(9)}$ & $\Q_{2}^{(3)}$ & $\Q_{p}^{(10)}$ & $\Q_{p}^{(9)}$ & \multicolumn{2}{|c|}{$\Q_{2}^{(3)}$} & $\Q_{p}^{(10)}$ & $\Q_{2}^{(3)}$ \\
\hline
$-2p(p-2)$ & $\Q_{2}^{(4)}$ & $\Q_{p}^{(9)}$ & \multicolumn{2}{|c|}{$\Q_{2}^{(4)}$} & $\Q_{p}^{(10)}$ & $\Q_{p}^{(9)}$ & $\Q_{2}^{(4)}$ & $\Q_{p}^{(10)}$ \\ 
\hline
\end{tabular}

\bigskip

\underline{Notes For Table}
\end{center}

\begin{enumerate}
\item If $b_{1} < 0$ and $b_{2} > 0$, equation $\eqref{eq1}$ has no solutions in $\mathbb{R}$.
\item If $b_{1} < 0$ and $b_{2} < 0$, equation $\eqref{eq2}$ has no sution in $\mathbb{R}$.
\item Suppose there exists a solution $(z_{1},z_{2},z_{3})$. We have $b_{1} \equiv 0$ (mod $2$) and $b_{2} \neq 0$ (mod $2$). Then comparing $2-$adic valuations of the left-hand and right-hand sides of equation $\eqref{eq1}$ easily implies $z_{1},z_{2} \in \Z_{2}$. But then $b_{1}z_{1}^2 - b_{1}b_{2}z_{3}^2 \equiv 0$ (mod $2$), so equation $\eqref{eq2}$ implies $p \equiv 0$ (mod $2$). Since $p$ is odd, this is a contradiction, so equations $\eqref{eq1}$ and $\eqref{eq2}$ have no solutions over $\Q_{2}$
\item Adding the $\Q_{2}$-non-trivial pairs from $(3)$ to the (pairs corresponding to the) points in $E[2]$ yields these $\Q_{2}$-non-trivial pairs.
\item If $(b_{1},b_{2})=(p,1)$ or $(p(p-2),1)$, then valuation arguments again show that any solution $(z_{1},z_{2},z_{3})$ has $z_{1},z_{2} \in \Z_{p}$, hence equation $\eqref{eq1}$ becomes 
\[
z_{2}^2 \equiv -2 \text{  (mod $p$)}
\]
which has no solutions since $p \equiv 7$ (mod $8$).
\item Adding the pairs from $(5)$ to the points in $E[2]$ yields these $\Q_{p}$-non-trivial points.
\item If $(b_{1},b_{2})=(p-2,1)$, then again $z_{1},z_{2} \in \Z_{p-2}$, so equation $\eqref{eq1}$ implies
\[
z_{2}^2 \equiv -2 \text{  (mod $(p-2)$),}
\]
which has no solutions since $p-2 \equiv 5$ (mod $8$).
\item Adding the pairs from $(7)$ yields these $\Q_{p-2}$-non-trivial points.
\item Suppose $b_{2} \equiv 0$ (mod $p$) and $b_{1} \neq 0$ (mod $p$), and suppose there exists a solution $(z_{1},z_{2},z_{3})$. Let
\[
k=\nu_{p}(z_{1}), \quad j=\nu_{p}(z_{2}, \quad l=\nu_{p}(z_{3}.
\]
Then equation $\eqref{eq2}$ implies that 
\begin{align*}
1 &= \nu_{p}(b_{1}z_{1}^{2} - b_{1}b_{2}z_{3}^2) \\
&= min\{2k, 1+2l\}
\end{align*}
which implies $l=0$ and $k > 0$. But equation $\eqref{eq1}$ implies 
\begin{align*}
0 &= \nu_{p}(b_{1}z_{1}^{2} - b_{2}z_{2}^2) \\
&=min\{ 2k,1+2j \}
\end{align*}
which implies $k=0$, which is a contradiction. Hence equations $\eqref{eq1}$ and $\eqref{eq2}$ have no solutions over $\Q_{p}$ in his case.

\item Adding pairs from (9) to the points in $E[2]$ yields these $\Q_{p}$-non-trivial pairs.

\item If $(b_{1},b_{2}) = (1,-(p-2))$, then once again equation $\eqref{eq1}$ implies $z_{1},z_{2} \in \Z_{p-2}$. Subtracting $p$ from both sides of equation $\eqref{eq1}$ yields
\begin{align*}
z_{1}^2 + (p-2)z_{2}^2 - p = 2-p,
\end{align*}
so looking modulo $(p-2)$ we have
\[
z_{1}^2 \equiv p \equiv 2 \text{  (mod $(p-2))$)}.
\]
But $p-2 \equiv 5$ (mod $8$), so no such $z_{1}$ exists, hence equation $\eqref{eq1}$ has no solutions over $\Q_{p-2}$.

\item Adding pairs from $(11)$ to the points in $E[2]$ yields these $\Q_{p-2}$-non-trivial pairs.

\item If $(b_{1},b_{2}) = (p,-(p-2))$ or $(p(p-2),-(p-2)$, then equation $\eqref{1}$ again implies $z_{1},z_{2} \in \Z_{p}$ and reduces to
\[
z_{2}^{2} \equiv -1 \text{  (mod $p$)}
\]
which has no solutions in $\Q_{p}$ since $p \equiv 3$ (mod $4$).

\item Adding the pairs from $(13)$ to the points in $E[2]$ yields these $\Q_{p}$non-trivial pairs.

\item If $(b_{1},b_{2}) = (p(p-2),-1)$, then equation $\eqref{eq1}$ implies $z_{1},z_{2} \in \Z_{p}$ and reduces to
\[
z_{2}^{2} \equiv 2 \text{  (mod $p-2$)}
\]
which has no solutions in $\Q_{p-2}$ since $p-2 \equiv 5$ (mod $8$).

\item Adding the pair from $(15)$ to the points in $E[2]$ yields these $\Q_{p-2}$-non-trivial pairs.
\end{enumerate}

\bigskip

This table reveals that the only elements of $\Q(S,2) \times \Q(S,2)$ in the image of the map $\phi$ are those we got from $E[2]$, and Theorem $2$ then implies
\[
E(\Q) / 2E(\Q) \hookrightarrow \Q(S,2) \times \Q(S,2) \cong (\Z / 2\Z)^{2}.
\]
Since we showed before that $E_{\text{tors}} = E[2]$, we know
\[
E(\Q) \cong \Z^{r} \times (\Z / 2\Z)^{2},
\]
hence
\[
E(\Q) / 2E(\Q) \equiv (\Z / 2\Z)^{2 + r}.
\]
Thus, we have shown that the rank of $E$ is $r=0$, and we have $E \cong (\Z / 2\Z)^{2}$.  $\quad  \square$

\bigskip

\section{Proving Theorem 1(b)}

Proving Theorem 1(a) was somewhat long and tedious, but it was not very difficult. To establish the bound of the rank of our curve when $p \equiv 5$ (mod $8$), we will need a more general form of $2-$descent. First we restate the theorem to be proved.

\bigskip

\textbf{Theorem 1(b).} \textit{Let $p$ and $p-2$ be prime numbers in $\Z$ with $p \equiv 5$ (mod $8$). Then the elliptic curve $E(\Q)$ given by}
\[
E: y^2 = x(x-p)(x-2)
\] 
\textit{has rank at most $1$. In particular, it is either of the form }
\[
E(\Q) \cong \Z / 2\Z \times \Z / 2\Z.
 \]
\text{or}
\[
E(\Q) \cong \Z \times \Z / 2\Z \times \Z / 2\Z.
 \]

\bigskip

To prove this, we use the method of $2$-descent described in chapter X of Silverman [1] as Proposition 4.9, which we now state.

\bigskip

\textbf{Theorem 3. (Descent via Two-Isogeny.)} \textit{Let $E/ \Q$ and $E' / \Q$ be elliptic curves given by equations}
\[
E: y^2 = x^3 + ax^2 + bx \text{   and   } E': Y^2 = X^3 - 2aX^2 + (a^2-4b)X;
\]
\textit{and let }
\[
\phi : E \to E' \quad \phi(x,y) = (y^2/x^2,y(b-x^2)/x^2)
\]
\textit{be the isogeny of degree $2$ with kernel $E[\phi]=\{{\cal{O}},(0,0)\}$. Let $S$ consist of $\infty$ and the places dividing $2b(a^2 - 4b)$, an let $\Q(S,2)$ be defined as before. There is an exact sequence}
\[
0 \rightarrow E'(\Q)/\phi(E(\Q)) \rightarrow \Q(S,2) \to WC(E/\Q)[\phi]
\]
\[
{\cal{O}} \to 1 \quad
\]
\[
(0,0) \to a^2 - 4b \quad d \to \{C_{d}/\Q\},
\]
\[
(X,Y) \to X \quad
\]
\textit{where $C_{d}/\Q$ is the homogeneous space for $E/ \Q$ given by the equation}
\[
C_{d} : dq^2 = d^2 - 2adz^2 + (a^2 - 4b)z^4.
\]
\textit{The $\phi-$Selmer group is then}
\[
S^{(\phi)}(E / \Q) \cong \{ d \in \Q(S,2) : C_{d}(\Q_{\nu}) \neq \varnothing \textit{ for all } \nu \in S\}.
\]

\bigskip

This method of $2-$descent is obviously more complicated than the first version we used. The theory behind this method is high-powered, and the interested reader is encouraged to see Silverman [1]. We will say only a few words about this theorem before using it to prove Theorem 1(b).

First, we did not previously define $WC(E / \Q)$; this is the Weil-Ch$\hat{a}$telet group for $E/ \Q$, which is the set of equivalence classes of homogeneous spaces for $E / \Q$. A homogeneous space for $E / \Q$ is a "twist" of the elliptic curve $E$, or some smooth curve on which $E$ has a simply transitive algebraic group action. What is important about these homogeneous spaces is that they encode certain information about the elliptic curve $E$, and Theorem 3 says that we are able to calculate the Selmer group $S^{(\phi)}$ of our isogeny $\phi$ by looking at these homogeneous spaces over local fields. 

What, then, is the Selmer group? This, too, should be investigated in Silverman [1]. Essentially, the Selmer group contains the homogeneous spaces which have $\Q_{\nu}$-rational points for every $\nu$.  The important thing about $S^{(\phi)}$ is that we have an exact sequence
\[
0 \rightarrow E'(\Q) / \phi(E(\Q)) \rightarrow S^{(\phi)}(E / \Q) \rightarrow \Sha (E / \Q ) [\phi] \rightarrow 0.
\]
The Shafarevich-Tate group $\Sha$ is essentially the group of homogeneous spaces for $E$ which have a $\Q_{\nu}$-rational point for every place $\nu$ but no $\Q-$rational points; again, see Silverman [1]. Together, the Selmer and Sha groups measure the failure of the Hasse principle for these curves. Finally, computing both $S^{(\phi)}$ and $S^{(\hat{\phi})}$ and using the above exact sequence together with the exact sequence
\[
0 \rightarrow \frac{E'(\Q)[\hat{\phi}]}{\phi(E(\Q)[2]} \rightarrow \frac{E'(\Q)}{\phi(E(\Q))} \rightarrow \frac{E(\Q)}{2E(\Q)} \rightarrow \frac{E(\Q)}{\hat{\phi}(E'(\Q))} \rightarrow 0 ,
\]
we can compute $E(\Q) / 2E(\Q)$ as before and thus deduce the rank of $E$. 

Before proving Theorem 1(b), we first prove two lemmas which will greatly simplify the proof of the theorem.

\bigskip

\textbf{Lemma 1. } \textit{Given the curves}
\[
E: y^2 = x^3 - (p+2)x^2 + 2px  \text{ and } E': y^2 = x^3 + 2(p+2)x^2 + (p-2)^2x
\]
\textit{and the isogeny}
\[
\phi: E \to E' \quad \phi(x,y) = \left( \frac{y^2}{x^2}, \frac{y(2p-x^2)}{x^2} \right),
\]
\textit{the selmer group $S^{(\phi)}$ is }
\[
S^{(\phi)} = \{1,-1\}.
\]

\bigskip

\textbf{Proof of Lemma 1. }
We first note that Ker $\phi = \{{\cal{O}},(0,0) \}$ and
\[
E[2] = \{ {\cal{O}}, (0,0),(2,0),(p,0) \} \textit{   and   } E'[2]=\{ {\cal{O}},(0,0) \}.
\]
Also, we have $2b(a^2-4b)=4p(p-2)^2$, so the set $\Q(S,2)$ is 
\[
\Q(S,2) = \{  \pm 1, \pm 2, \pm p, \pm (p-2), \pm 2p, \pm 2(p-2), \pm p(p-2), \pm 2p(p-2) \}.
\]
For each $d \in \Q(S,2)$, we must check whether the associated homogeneous space
\[
C_{d}: dw^2 = d^2 + 2(p+2)dz^2 + (p-2)^2 z^4
\]
has points over each local field $\Q_{q}$ for $q=2,p,p-2$. From Theorem 3, we have
\[
{\cal{O}}, (0,0) \to 1
\]
and so $d=1 \in$ $S^{(\phi)}$, and it remains to check whether the remaining $d \in \Q(S,2)$ are in the Selmer group. 

\bigskip
\begin{enumerate}
\item $d= \pm p$ 

We present the argument for $C_{p}$ to show that $d=p$ is not in the Selmer group; the argument for $C_{-p}$ is identical. Our homogeneous space is
\[
C_{p} : pw^2 = p^2 + 2(p+2)pz^2 + (p-2)z^4. 
\]
Now, we have that the $p-$adic valuation of the left-hand side, $\nu_{p}$(LHS), is odd. Hence, the valuation of the right-hand side must also odd. But we have
\[
\nu_{p}\text{(RHS)} = \text{min}\{2, 1+2k, 4k \} \text{    where } k=\nu_{p}(z).
\]
(We have equality in the above expression since no two of $2,1+2k,$ and $4$ can ever be the same.) If $k \leq 0$, then this minimum equals $4k$, which is even. Similary, if $k>0$, then the minimum is $2$. In neither case can the valuations of the left- and right-hand sides match, so $C_{p}$ has no $\Q_{p}$-rational points. Hence $\pm p \notin S^{(\phi)} $.

\item $d = \pm 2$ 

We present the argument for $C_{2}$ to show that $d=2$ is not in the Selmer group; the argument for $C_{-2}$ is identiical. Our homogeneous space is
\[
C_{2}: 2w^2 = 2^2 - 2^2(p+2)z^2 + (p-2)^2z^4.
\]
The $2-$adic valuation of the left hand side is odd, while
\[
\nu_{2}\text{(RHS)} \geq \text{min} \{2, 2+2k, 4k \} \text{   where } k=\nu_{2}(z).
\]
Since $\nu_{2}$(RHS) must be odd, we must have at least two of $2$, $2+2k$, and $4k$ be equal and minimal; this is easily seen to be impossible, hence $C_{2}$ has no $\Q_{2}$-rational points. Thus $\pm 2 \notin S^{(\phi)}$.

\item $d= \pm 2p$ 

Either of the two arguments above show that $\pm 2p \notin S^{(\phi)}$.

\bigskip

\noindent Notice that since $S^{(\phi)}$ is a group, we have
\[
\pm 2(p-2), \pm 2p(p-2), \pm p(p-2) \notin S^{(\phi)}.
\]
\noindent It remains to check whether $C_{d}$ has $\Q_{q}$-rational points for 
\\ $d \in \{ -1, \pm (p-2) \}$ and $q\in \{2,p,p-2\}$.

\item $d=-1$ 

Our homogeneous space is
\[
C_{-1} : -w^2 = 1-2(p+2)z^2 + (p-2)^{2}z^{4}.
\]
We show that $C_{-1}$ has $\Q_{q}$-rational points for each $q \in \{2, p, p-2\}$, hence $-1 \in S^{(\phi)}$.

\begin{enumerate}
\item $q=2$
We find a solution in $\Q_2$ as follows. Let 
\[
f(w,z)=1-2(p+2)z^2 + (p-2)^{2}z^{4} + w^{2}.
\]
Then we have
\begin{align*}
f(2,1) &= 1 - 2(p+2) + (p-2)^{2} + 4 \\
&= 5 - 2(7+8k) + (3+8k)^2 \\
&\equiv 5-14-16k+9+16k \text{   (mod 32)  since } p \equiv 5 \text{  (mod $8$)} \\
& \equiv 0 \text{ (mod $32$)}.
\end{align*}
So $\nu_{2}(f(2,1)) \geq 5$, but 
\[
\frac{\partial f}{\partial w} = 2w \bigg|_{(w,z)=(2,1)} = 4,
\]
so $2 \nu_{2}(\partial f / \partial w)=4 < 5$. Hence by Hensel's lemma this lifts to a solution in $\Q_{2}$.

\item $q=p$

We find a solution in $\Q_p$ as follows. With $f(w,z)$ as before, we have
\[
f(w,z) = 1 - 4z^2 + 4z^4 + w^2 \text{  (mod $p$)},
\]
hence
\[
f(w_{0},0) = 1 + w^2 \equiv 0 \text{  (mod $p$)}
\]
has a solution $w_{0}$ since $p \equiv 4$ (mod $8$). Since $p \not| w_{0}$ and $p$ is an odd prime, we have
\[
\frac{\partial f}{\partial w}=2w \bigg|_{(w,z)=(w_{0},0)} \neq 0 \text{  (mod $p$)},
\]
so by Hensel's lemma this lifts to a solution in $\Q_{p}$.

\item $q=p-2$

We find a solution in $\Q_{p-2}$ as follows. With $f(w,z)$ as before, we have
\[
f(w,z) \equiv 1 - 2(p+2)z^2 + w^{2} \text{  (mod $(p-2)^2$)}
\]
so
\[
f(w,0) \equiv 1 + w^2 \equiv 0 \text{  (mod $(p-2)^2$)}
\]
has a solution $w_{0}$ since $(p-2)^2 \equiv 1$ (mod $4$). Once again, we have
\[
\frac{\partial f}{\partial w} = 2w \bigg|_{(w,z)=(w_{0},0)} \neq 0 \text{   (mod $(p-2)$)}
\]
so by Hensel's lemma this lifts to a solution in $\Q_{p-2}$.

\end{enumerate}

\bigskip

\noindent Thus, we have $-1 \in S^{(\phi)}$. We now show that $C_{-(p-2)}$ has no $\Q_{2}-$rational points, hence $-(p-2) \notin S^{(\phi)}$. Since $S^{(\phi)}$ is a group, this will imply $(p-2) \notin S^{(\phi)}$, completing our calculation of $S^{(\phi)}$. 

\item $d = -(p-2)$

Our homogeneous space is
\[
C_{-(p-2)} : -(p-2)w^{2} = (p-2)^{2} - 2(p+2)(p-2)z^{2} + (p-2)z^{4}.
\]
We will show this space has no points over $\Q_{2}$. Let $\nu_{2}(w)=k, \nu_{2}(z)=j$. Suppose there is a solution $(w,z)$. We have
\[
\nu_{2}\text{(LHS)} = 2k
\]
and
\[
\nu_{2}\text{(RHS)} \geq \text{min} \{0, 1+2j, 4j \}.
\]
There are two cases to consider.

Suppose $j>0$, which implies $k=0$. Then $w,z \in \Z_{2}$ with $w$ odd and $z$ even, so
\[
f(w,z) \equiv (p-2)^{2} + (p-2)w^{2} \equiv 1+ 3w^2 \text{  (mod $8$)}
\]
which has no solutions. Having no solutions modulo $8$ implies that there are no solutions in $\Q_{2}$.

Now suppose $j=0$, which implies $k \geq 0$. Then $w,z \in \Z_{2}$ and $z$ is odd. We examine $f(w,z)$ modulo different powers of $2$ to gain information about $w$ and $z$. Remember that $p \equiv 5$ (mod $8$). Looking modulo $2$ yields
\[
f(w,z) \equiv 1 + 1 + w^2 \quad \text{(mod $2$)}
\]
since $z$ is odd, so this implies $w$ is even.

\bigskip

Looking modulo $2^2$ and $2^3$ yields no new information, so we look modulo $2^4$. Note that $(p-2)^2 \equiv 9$ and $(p+2)(p-2) \equiv 5$ modulo $2^4$. Then, writing $z=1 + 2r$, $w=2s$, and $p = 5 + 8l$, we have
\begin{align*}
f(w,z) &\equiv 9 - 10(1+2r)^2 + 9(1+2r)^4 + (3+8l)(2s)^2 \quad \text{(mod $16$)} \\
& \equiv 8 + 12s^2 \quad \text{(mod $16$)}
\end{align*}
and having this equal to $0$ modulo $16$ is equivalent to having
\[
2 + 3s^2 \equiv 0 \quad \text{(mod $4$)}
\]
which has no solutions. This implies that there are no solutions to $f(w,z)=0$ in $\Q_{2}$.

Now suppose $j<0$, which implies $k=2j$. Changing the signs of $j$ and $k$, write $w=2^{-2j}w_{0}$ and $z=2^{-j}z_{0}$, where $w_{0}$ and $z_{0}$ are both odd. The equation defining our homogeneous space is then
\[
-(p-2)2^{-4j}w_{0}^{2} = (p-2)^{2} - 2(p+2)(p-2)2^{-2j}z_{0}^2 + (p-2)^{2}2^{-4j}z_{0}^{4}.
\]
Multiplying through by $2^{4j}$, this becomes
\[
-(p-2)w_{0}^2 = 2^{4j}(p-2)^{2} - (p+2)(p-2)2^{2j+1}z_{0}^{2} + (p-2)^{2}z_{0}^{4}, 
\]
and looking at the resulting function modulo $8$ yields
\begin{align*}
f(w,z) &\equiv (p-2)^{2}z_{0}^{4} + (p-2)w_{0}^{2} \quad \text{(mod $8$)} \\
&\equiv z_{0}^{4} + 3w_{0}^{2} \quad \text{(mod $8$)}.
\end{align*}
Since $z_{0},w_{0}$ are odd, we have $f(w,z) \equiv 4 \neq 0$ modulo $8$. Since there are no solutions modulo $8$, there are no solutions in $\Q_{2}$.

\end{enumerate}

We conclude that $-(p-2) \notin S^{(\phi)}$, and since $-1 \in S^{(\phi)}$, we have $(p-2) \notin S^{(\phi)}$. Thus, $S^{(\phi)} = \{ 1, -1 \}.$ $\quad \square$

\bigskip

\textbf{Lemma 2. }
\textit{Given the curves}
\[
E':Y^2 = X^3 + 2(p+2)X^2 + (p-2)X \textit{ and } E:y^2 = x^3 - (p-2)x^2 + 2px
\]
\textit{and the isogeny}
\[
\hat{\phi}: E' \to E \quad \hat{\phi}(X,Y) = \left( \frac{y^2}{x^2}, y \frac{p-2-x^2}{x^2}\right)
\]
\textit{the selmer group Sel$^{(\hat{(\phi)})}$ is }
\[
S^{(\hat{\phi})} = \{ 1,2,p,2p\}.
\]

\bigskip

\textbf{Proof of Lemma 2.}

We first note that Ker$\hat{\phi}=\{{\cal{O}},(0,0)\}$ and
\[
E[2] = \{ {\cal{O}}, (0,0),(2,0),(p,0) \} \text{   and   } E'[2]=\{ {\cal{O}},(0,0) \}.
\]
Our set $\Q(S,2)$ is once again
\[
\Q(S,2) = \{  \pm 1, \pm 2, \pm p, \pm (p-2), \pm 2p, \pm 2(p-2), \pm p(p-2), \pm 2p(p-2) \}.
\]
For each $d \in \Q(S,2)$, we must check whether the associated homogeneous space
\[
C_{d}: dw^2 = d^2 - 4(p+2)dz^2 + 2pz^4
\]
has points over each local field $\Q_{q}$ for $q=2,p,p-2$. From Theorem 3, we have
\[
{\cal{O}} \to 1 \quad (0,0) \to 2p \quad (2,0) \to 2 \quad (p,0) \to p
\]
and so $\{1,2,p,2p\} \subset S^{\hat{\phi}}$; it remains to check whether the remaining $d \in \Q(S,2)$ are in the Selmer group. The set of $d$ which remain to be checked is
\[
\{ -1, -2, -p, \pm (p-2), \pm 2p, \pm 2(p-2), \pm p(p-2), \pm2p(p-2) \}.
\]
We proceed as in the last lemma, only this time less work is necessary. Consider $d= \pm (p-2)$; we show that $C_{p-2}$ has no solutions over $\Q_{p-2}$, and the argument for $d=-(p-2)$ is identical. Our homogeneous space is
\[
C_{p-2} : \pm (p-2)w^2 = (p-2)^2 \pm 4(p+2)(p-2)z^2 + 2pz^4
\]
Letting $k = \nu_{p-2}(w) = k$ and $\nu_{p-2}(z)=j$, we have $\nu_{2}$(LHS)$=1+2k$ and \\ 
$\nu_{2}$(RHS)$ \geq$ min$\{2,1+2j,4j\}$, which implies $j=k=0$, hence $w,z \in \Z_{p-2}$, and in particular, $(p-2)$ does not divide $w,z$. Taking as our function
\[
f(w,z) = (p-2)^2 \pm 4(p+2)(p-2)z^{2} + 2pz^4 \mp (p-2)w^2,
\]
 then looking at the function modulo $p-2$ yields 
\begin{align*}
f(w,z) &\equiv 4z^4 \quad \text{(mod $(p-2)$)} 
\end{align*}
which has no solutions $(w,z)$ with $p-2 \not| z$, hence $f(w,z)$ has no solutions in $\Q_{\pm (p-2)}$. We conclude that $\pm(p-2) \notin S^{(\hat{\phi})}$. But by the group structure of $S^{(\hat{\phi})}$, this eliminates all of the other possible $d \in \Q(S,2)$. Hence we have $S^{(\hat{\phi})} = \{ 1,2,p,2p\}$. $\quad \square$

\bigskip

We now prove the theorem.

\bigskip

\textbf{Proof of Theorem 1(b).} 

First we show that $E_{\text{tors}} = E[2] \cong (\Z / 2Z)^2$. Suppose $p > 5$. As in the proof of Theorem 1(a), $3$ is a prime of good reduction. Since $p \equiv 5$ (mod $8$), we have $p \equiv 2 + 2k$ (mod $3$) for some $0 \leq k leq 2$. But if $k=2$, then $p \equiv 0$ (mod $3$) and is not prime (since $p > 5$). Similarly, if $k=0$, then $p -2 \equiv 0$ (mod $3$) and so is not prime (since $p >5$ implies $p-2 > 3$). Then one easily checks that
\[
E(\mathbb{F}_3) = \{  {\cal{O}}, (0,0), (1,0), (2,0) \},
\] 
and since $E_{{tors}} \hookrightarrow E(\mathbb{F}_3)$ and $E[2] \subset E_{\text{tors}}$, we have 
\[
E_{\text{tors}} = E[2] \cong (\Z / 2Z)^2.
\]
If $p=5$, one easily checks that
\[
E(\mathbb{F}_7) = \{ {\cal{O}}, (0,0), (1,2), (1,5), (2,0), (3,1), (3,6), (5,0)  \} \cong (\Z / 4\Z) \times (\Z / 2\Z)
\]
and a similar argument gives the desired result.

Now we compute the rank. To apply Theorem 3, we need to compute precisely the Selmer groups from Lemmas 1 and 2. Having computed both Selmer groups, we may now use the exact sequences mentioned earlier to compute the rank of $E$. Recall that we have the following exact sequences:

\bigskip

\begin{equation}\label{ex1}
0 \to \frac{E'}{\phi(E)} \to S^{\phi}(E) \to \Sha (E)[\phi] \to 0
\end{equation}
\begin{equation}\label{ex2}
0 \to \frac{E}{\hat{\phi}(E)} \to S^{\hat{\phi}}(E) \to \Sha (E')[\hat{\phi}] \to 0
\end{equation}
\begin{equation}\label{main}
0 \to \frac{E'[\hat{\phi}]}{\phi(E[2])} \to \frac{E'}{\phi(E)} \to \frac{E}{2E} \to \frac{E}{\hat{\phi}(E')} \to 0
\end{equation}

\bigskip

Now, we have computed
\[
S^{(\phi)} \cong \Z / 2\Z \quad \text{ and } \quad S^{(\hat{\phi})} \cong (\Z / 2\Z)^2,
\]
and we actually showed that every point in $S^{(\hat{\phi})}$ came from a point on our curve, hence $\Sha (E')[\hat{\phi}] = 0$. Sequence $\eqref{ex2}$ then implies

\[
E / \hat{\phi}(E) \cong (\Z / 2\Z)^2.
\]

It is also easy to check that $\phi(E[2])=E'[\hat{\phi}]$, hence $E'[\hat{\phi}] / \phi(E[2]) = 0$. Thus, sequence $\eqref{main}$ becomes

\begin{equation}\label{main2}
0 \to \frac{E'}{\phi(E)} \to \frac{E}{2E} \to (\Z / 2\Z)^2 \to 0.
\end{equation}

We do not know if $\Sha(E)[\phi] = 0$; thus, from sequence $\eqref{ex1}$ we have

\[
0 \to \frac{E'}{\phi(E)} \to \Z / 2\Z \to \Sha (E)[\phi] \to 0,
\]

\noindent so we have either $\frac{E'}{\phi(E)} \cong \Z / 2\Z$ or $0$; comparing the orders of the groups in $\eqref{main2}$ shows that rank$(E) = 1$ in the first case and $0$ in the second, proving the theorem. $\quad \square$

\section{Acknowledgements}
The author would like to thank Dr. Thomas Hagedorn for all of his help and guidance with this paper, and Jason Beers for making the conjectures that feature as this paper's theorems. This paper was written in partial fulfilment of graduation requirements at The College of New Jersey.

\section{References}

[1] Joseph H. Silverman. \textit{The Arithmetic of Elliptic Curves.} Springer-Verlag. New York, NY, 1986.
\end{document}